\documentclass[11pt]{article}
\usepackage{latexsym,amssymb}
\usepackage{enumerate}
\usepackage{amsthm}
\usepackage{amsmath}
\newtheorem{corollary}{Corollary}[section]

\newtheorem{lemma}[corollary]{Lemma}
\newtheorem{proposition}[corollary]{Proposition}
\newtheorem{remark}[corollary]{Remark}

\newtheorem{theorem}[corollary]{Theorem}
\newcommand{\mylabel}[1]{\label{#1}
            \ifx\undefined\stillediting
            \else \fbox{$#1$}\fi }
\newcommand{\BE}{\begin{equation}}
\newcommand{\BEQ}[1]{\BE\mylabel{#1}}
\newcommand{\EEQ}{\end{equation}}
\newcommand{\rfb}[1]{\mbox{\rm
   (\ref{#1})}\ifx\undefined\stillediting\else:\fbox{$#1$}\fi}

\newfont{\Blackboard}{msbm10 scaled 1200}

\newfont{\roma}{cmr10 scaled 1200}
\def\cal{\mathcal}

\def\cA{{\cal A}}
\def\cB{{\cal B}}

\def\cD{{\cal D}}

\def\cH{{\cal H}}

\newcommand{\half}   {{\frac{1}{2}}}
\newcommand{\Dscr}   {{{\mathcal D}}}
\newcommand{\mm}    {{\hbox{\hskip 0.5pt}}}
\newcommand{\m}     {{\hbox{\hskip 1pt}}}

\newcommand{\bluff} {{\hbox{\raise 15pt \hbox{\mm}}}}
\newcommand{\R}   {\mathbb{R}}

 \allowdisplaybreaks \textwidth=16cm \textheight=23cm
\makeatletter
\def\section{\@startsection {section}{1}{\z@}{-3.5ex plus -1ex minus
    -.2ex}{2.3ex plus .2ex}{\large\bf}}
\makeatother
\def\be{\begin{equation}}
\def\ee{\end{equation}}

\begin{document}

\title
{\bf Stabilization of coupled second order systems with delay}
\footnotetext[1]{Research supported by the IRD-UMMISCO,  and by CMPTM-Research project : 10/MT /37.}

\author{E. M.  Ait Benhassi $^\ddag$, K.  Ammari
\thanks{D\'epartement de Math\'ematiques, Facult\'e des
Sciences de Monastir, Universit\'e de Monastir, 5019 Monastir,
Tunisie, e-mail: kais.ammari@fsm.rnu.tn},  S. Boulite \thanks{D\'epartement de Math\'ematiques, Facult\'e des
Sciences Ain Chock, Universit\'e   Hassan II,  Casablanca, e-mail s.boulite@fsac.ac.ma}   and
L. Maniar
\thanks{Universit\'e Cadi
Ayyad,  Facult\'e des Sciences Semlalia, LMDP, UMMISCO (IRD- UPMC), 
Marrakech 40000, B.P. 2390, Maroc, Fax. 0021224437409, e-mail:
m.benhassi@ucam.ac.ma, 
maniar@ucam.ac.ma}}
\date{}

\maketitle

\begin{abstract}
In this paper we characterize the output feedback stabilization of some
 coupled systems with delay. The proof of the main result uses the method
introduced in Ammari and Tucsnak \cite{at} where the exponential
stability for the closed loop system is reduced to an observability
estimate for the corresponding conservative adjoint  system, under  a
boundedness condition of the transfer function of the associated open
loop system.
\end{abstract}

\section{Introduction}

In this paper, our purpose is to characterize the output feedback  stabilization of
coupled second order  infinite dimensional systems by only one 
feedback. Using an output  feedback,  the closed loop system we treat is the  following
\be\label{cwo1} \ddot w_1(t) + A_1w_1(t) + BB^*\dot w_1(t)+ C\dot  w_2(t) = 0, \quad t\geq 0, 
\ee
\be\label{cwo2} \ddot w_2(t) + A_2w_2(t) -C^*\dot w_1(t) = 0, \quad t\geq 0,\ee
\be \label{cwo3} w_i(0) =\m w_i^0, \,\dot w_i(0)=w_i^1,\,i=1,2. 
\ee
 Here, the operators $A_1,A_2$ are  unbounded positive self adjoint operators in Hilbert spaces $H_1,H_2$, respectively. 
The control operator $B$, acting only on the first equation,  is assumed here to be unbounded from 
$U$, another Hilbert space,  to $D(A_1^{\frac12})^*$.
The coupling operator $C$ is not necessarily bounded.
In \cite{Bous1,Bous2}, the authors have considered
coupled  systems in the case of bounded (even compact) 
coupling  operators $C$. In this  case  
the exponential stability does not  hold, since the  equation \eqref{cwo3} is  conservative when $C=0$. In stead, they studied the polynomial stability. 
Recently, Ammari and Nicaise \cite{Am-Ni} have characterized the
exponential energy decay  of these systems by an observability inequality of associated conservative adjoint systems,
 augmented 
with the 
output $y(t)=-B^*\dot w_1(t)$, in the case of bounded coupling  operators $C$. In \cite{Souf}, the author 
studied also these coupled systems in the case of unbounded coupling operators,  consedering  bounded operators $B$. 
 In this paper, we assume that both operators $B$ and $C$ are 
unbounded, and show the same result as in \cite{Am-Ni} using different arguments. Here, we transform the system  
\eqref{cwo1}-\eqref{cwo3} to a second order equation
\be \label{wavetaxi}
  \ddot w (t)+A w(t)+\mathcal{B}_0\mathcal{B}_0^*\dot w(t)=0,\quad t\geq0,
 \ee
  \be \label{wavetaxi1} 
 w(0)=w^0,\,\dot w(0)=w^1
\ee 
in the product space $H:=H_1\times H_2$, with appropriate operators $A$ and $\mathcal{B}_0$ defined in \eqref{a,b}. 
Then, we use the result of Ammari-Tucsnak \cite[Theorem 2.2]{at} to characterize the exponential enregy 
decay of  the  equation \eqref{wavetaxi}-\eqref{wavetaxi1}, 
and then deduce the one of  the coupled 
systems \eqref{cwo1}-\eqref{cwo3}.

The second aim of this paper is to characterize the exponential energy decay of 
the following coupled systems with delay
\be\label{del1} \ddot w_1(t) + A_1w_1(t) +  \alpha_1 \, BB^*\dot w_1(t)+ \alpha_2 \,
BB^*\dot w_1(t-\tau)+ C\dot  w_2(t) = 0, \quad t\geq 0, \ee
\be \ddot w_2(t) + A_2w_2(t) -C^*\dot w_1(t) = 0, \quad t\geq 0,\ee
\be \label{del3} w_i(0) = w_i^0, \,\dot w_i(0)=w_i^1,\,i=1,2,\,\, \dot w_1(s) = f_0(s), \,\,s\in
(-\tau,0). \ee
The operators $A_i, i=1,2$, $B,C$ satisfy the same conditions as above, and $\alpha_1,\alpha_2$ are positive constants. 
The introduction of a delay term in partial differential equations and its effect on the stabilization of these equations were the 
subjet of several papers, see for instance, \cite{abmb, batkai1,batkai, dat1,dat2,dat3,Nicaise1,Nic}, and the references therein. 
By the same technic as for the first coupled systems, we tranform  the system \eqref{del1}-\eqref{del3} to a second order equation with delay
\be \notag
  \ddot w (t)+A w(t)+\alpha_1\mathcal{B}_0\mathcal{B}_0^*\dot w(t)+
\alpha_2\mathcal{B}_0\mathcal{B}_0^*\dot w(t-\tau)=0,\quad t\geq0,
  \notag\ee
  \be \label{initialesb} w(0)=w^0,\,\dot w(0)=w^1, \dot w_1(s) = f_0(s), \,\,s\in
(-\tau,0).
\notag\ee 
At this level, our results in \cite{abmb} will allow us to conclude. 

We then apply our abstract results to two systems of coupled string  equations with delay. 
The first example is a coupled two string  
equations with ponctuel control and Dirichlet  boundary conditions
\begin{align*}
&\ddot w_1(t,x)-\frac{\partial^2 w_1}{\partial x^2}(t,x)+\alpha_1 \dot w_1(t,\xi)\delta_\xi+\alpha_2 \dot
 w_1(t-\tau,\xi)\delta_\xi+\beta \frac{\partial \dot w_2}{\partial x}(t,x)=0, (t,x)\in(0,\infty)\times (0,1),\\
&\ddot w_2(t,x)-\frac{\partial^2 w_2}{\partial x^2}(t,x)+\beta \frac{\partial \dot w_1}{\partial x}(t,x)=0,
\quad (t,x)\in(0,\infty)\times (0,1),\\
&w_i(t,0)=w_i(t,1)=0,\quad  t\in (0,\infty),\; i=1,2, \\
&w_i(0,x)=w_i^0(x),\,\,\dot w_i(0,x)=w_i^1(x),\,\,\dot w_1(s,x)=f_0(s,x), -\tau\leq s < 0,\, \; x\in (0,1), \,\;
i=1,2,\label{Ap1c}
\end{align*}
with $\xi\in (0,1)$, $\beta>0$ and $0<\alpha_2<\alpha_1$. We show that this system is not exponentially stable for all 
$\xi\in (0,1)$ and $\beta>0$, showing that the observability  inequality of its conservative adjoint  system can not hold.  
To give a positive application of our abstract results, we  consider a coupled two wave 
equations with ponctuel control and mixed  boundary conditions
\begin{align*}
&\ddot w_1(t,x)-\frac{\partial^2 w_1}{\partial x^2}(t,x) + w_1(t,x) + \alpha_1 \dot w_1(t,\xi)\delta_\xi+\alpha_2 \dot
 w_1(t-\tau,\xi)\delta_\xi+\beta \frac{\partial \dot w_2}{\partial x}(t,x)=0,\;  t\geq 0,x\in (0,1),\\
&\ddot w_2(t,x)-\frac{\partial^2 w_2}{\partial x^2}(t,x) + w_2(t,x) + \beta \frac{\partial \dot w_1}{\partial x}(t,x)=0,\,
 t\geq 0,x\in (0,1),\\
&\frac{\partial w_1}{\partial x} (t,0) = \frac{\partial w_1}{\partial x}(t,1) = 0, \, w_2(t,0)=w_2(t,1)=0,\, t\geq 0, \\
&w_i(0,x)=w_i^0(x),\,\,\dot w_i(0,x)=w_i^1(x),\,\dot w_1(s,x)=f_0(s,x), -\tau\leq s < 0,x\in (0,1),\;  i=1,2
\end{align*}
with $\xi\in (0,1)$, $\beta$ is a positive constant and  $0<\alpha_2<\alpha_1$. Using the 
classical inequality by Ingham \cite{Ing} for non-harmonic Fourier series, 
we show that the 
observability  inequality of the conservative adjoint  system holds if and only if 
$\xi$ is a rational number with coprime factorisation $\xi = \frac{p}{q}$, where $p$ is odd. Thus, 
this is a necessary and sufficient 
condition for the  exponential energy decay  of the above system.

\section{Problem formulation }

Let  $H_i$  be a Hilbert space equipped with the norm
$||\cdot||_{H_i},$\,$i=1,2$ and let \be\label{as00} A_i~:~H_i~\supseteq~ \cD(A_i)\rightarrow~H_i, i=1,2,  \mbox{ be \,positive\, self\, adjoint\, operators.}\ee  We introduce the scale of Hilbert spaces $H_{i,{\alpha}}$ as
$H_{i,{\alpha}}=\Dscr(A_i^{\alpha})$  with the norm $\|z
\|_{i,{\alpha}}=\|A_i^\alpha z\|_{H_i}$ and  their dual spaces $H_{i,{-\alpha}}=H_{i,{\alpha}}^*,$\,$i=~1,2$.
The second ingredient needed for our construction is a control operator $B$ such that \be\label{as0} B~ :~ U
\longrightarrow ~H_{1,-\half}\,\,\mbox{is\, bounded,  }\,\ee where $U$ is another Hilbert space
identified with its dual. The operator $B^*$ is then bounded from
$H_{1,\frac12}$ to $U$.
 We need also  a unbounded
linear operator $C: H_2 \supseteq~ \cD(C)\longrightarrow ~H_{1}$  satisfying   the following assumptions
  \be\label{as1}H_{1,\half}\hookrightarrow \cD(C^*) \,\, \mbox{and}\,\,\,H_{2,\half}\hookrightarrow \cD(C).\ee
\begin{remark}\label{remarque}
By assumptions \eqref{as1}, one can see that the operators  $CA_2^{-\half}$ and $ C^*A_1^{-\half}$ can be extended to bounded operators from $H_2$ to $H_1.$
\end{remark}

The first coupled systems that  we consider are described by \BEQ{damped01} \ddot
w_1(t) + A_1w_1(t) +   BB^*\dot w_1(t)+ C\dot  w_2(t) = 0, \quad t\geq 0, \EEQ
\BEQ{damped02} \ddot w_2(t) + A_2w_2(t) -C^*\dot
w_1(t) = 0, \quad t\geq 0,\EEQ
 \BEQ{output0} w_i(0) \m=\m w_i^0,  \m
\dot w_i(0)=w_i^1,\,i=1,2, \EEQ where  the initial data $(w_1^0,w_1^1,w_2^0,w_2^1)$
belongs to a suitable space.

The equation \rfb{damped01} is understood as equation in 
$H_{1,{-\half}}$, i.e., all the terms are in $H_{1,{- \half}}$. The term
$BB^* \dot w_1(t)$ represents a feedback damping. 
Transforming system \eqref{damped01}-\eqref{output0} on  a second order system and using the method in \cite{at}, 
we  characterize the stabilization of this system. Namely, assuming that there exists $\delta\in[0,\half)$ such that for all $(x,y)\in H_{1,1}\times H_{2,1}$

\be\label{condition}|<x, C y>|\leq \delta\left(\|A_1^{\half}x\|^2+\|y\|^2+\|C^*x\|^2\right),\ee
  under the boundedness of corresponding  transfer function, system \eqref{damped01}-\eqref{output0} 
is exponentially stable if and only if there exists a constant $c>0$ such that \begin{align} \label{io02} c\int_0^T\|(B^*\phi)'\|_U^2dt\geq\|A_1^{\half}\phi(0)\|^2+\|A_2^{\half} \psi(0)\|^2+\|     \left(  \begin{smallmatrix}
                \dot \phi(0) \\
                \dot \psi(0) \\
              \end{smallmatrix}
            \right)\|^2_{H_1\times H_2}\end{align}

 for all  solution  $(\phi,\psi)$  of the following conservative adjoint  system
\be\notag \ddot \phi + A_1\phi +C\dot\psi=0\ee
\be\notag \ddot \psi+A_2 \psi-C^*\dot\phi=0.\ee

Our second interest  is  to characterize the stabilization of the following coupled systems with delay
\BEQ{damped1} \ddot
w_1(t) + A_1w_1(t) +  \alpha_1 \, BB^*\dot w_1(t)+ \alpha_2 \,
BB^*\dot w_1(t-\tau)+ C \dot w_2(t) = 0, \quad t\geq 0, \EEQ
\BEQ{damped2} \ddot w_2(t) + A_2w_2(t) -C^* \dot
w_1(t) = 0, \quad t\geq 0,\EEQ
 \BEQ{output} w_i(0) \m=\m w_i^0,  \m
\dot w_i(0)=w_i^1,\,i=1,2, \,\, \dot w_1(s) = f_0(s), \,\,s\in
(-\tau,0), \EEQ where $\tau
> 0$ is the time delay, $\alpha_1$ and $\alpha_2$ are positive real
numbers, and the initial data $(w_1^0,w_1^1,w_2^0,w_2^1,f_0)$
belongs to a suitable space. Assuming that $\alpha_2< \alpha_1$, under the same assumption \eqref{condition}  we prove that \eqref{damped1}-\eqref{output} is exponentially stable if and only if the observability inequality \eqref{io02} is satisfied, which is then equivalent to the exponential stability of \eqref{damped01}-\eqref{output0}.

\section{Coupled second order systems without delay}
Consider the following coupled systems
\BEQ{damped001} \ddot
w_1(t) + A_1w_1(t) +   BB^*\dot w_1(t)+ C\dot  w_2(t) = 0, \quad t\geq 0, \EEQ
\BEQ{damped002} \ddot w_2(t) + A_2w_2(t) -C^*\dot
w_1(t) = 0, \quad t\geq 0,\EEQ
\BEQ{output00} w_i(0) \m=\m w_i^0,  \m
\dot w_i(0)=w_i^1,\,i=1,2. \EEQ
 After  studying  the well-posedness of the coupled systems \eqref{damped001}-\eqref{output00}, we give a characterization of its exponential stability.
\subsection{Well-posedness}
Some change of variables,  leads to the following result
\begin{theorem}\label{transform}
If $(w_1,w_2)$ is a solution of \eqref{damped001}-\eqref{output00}, then  $(u,v)$ defined by
$$u=w_1,\;\; v=A_2^{-\half}\dot w_2-A_2^{-\half}C^*w_1$$ is a solution of
the system
 \BEQ{damped u}
 \ddot u(t) +(A_1+CC^*)u(t)+ BB^*\dot u(t)+ CA_2^{\half}v(t)=0 , \quad t\geq 0,
 \EEQ
 \BEQ{damped v}
 \ddot v(t)+ A_2 v(t) + A_2^{\half}C^*u(t)=0 , \quad t\geq 0.
 \EEQ
 \BEQ{outputuv}
  u(0)=u^0, \dot u(0)=u^1, v(0)=v^0, \dot v(0)=v^1
\EEQ
with
$
 u^0=w_1^0,  u^1=w_1^1, v^0=A_2^{-\half} w_2^1-A_2^{-\half}C^*w_1^0,  v^1=-A_2^{\half}w_2^0.
 $

Conversely, if $(u, v)$ is a solution of \eqref{damped u}-\eqref{outputuv}, then $(w_1,w_2)$ defined by
$$w_1=u,\;\; w_2=-A_2^{\half}\dot v$$
 is a solution of \eqref{damped001}-\eqref{output00}.
 \end{theorem}
 \proof
 Let $(w_1,w_2)$ be a solution of \eqref{damped001}-\eqref{output00}.
 Setting $u=w_1\,\,\mbox{ and}\,\,v=A_2^{-\half}\dot w_2-~A_2^{-\half}~C^*w_1,$ 
we have
 \be \notag u(t)=w_1(t),\;\; v(t)=A_2^{-\half}\dot w_2(t)-A_2^{-\half}C^*w_1(t),\quad t\geq 0,\ee
 \be \notag \dot u(t)=\dot w_1(t),
\;\; \dot v(t)=A_2^{-\half}\ddot w_2(t)-A_2^{-\half}C^*\dot w_1(t),\quad t\geq 0.\ee
 Equation \eqref{damped002} yields 
 \be \label{tr1} u(t)=w_1(t),\;\; v(t)=A_2^{-\half}\dot w_2(t)-A_2^{-\half}C^*w_1(t),\;\; \dot u(t)=\dot w_1(t),\quad t\geq 0,\,\dot v(t)=-A_2^{-\half}w_2(t),\quad t\geq 0.\notag\ee
 Thus
\be\label{change}
\left(\begin{smallmatrix}u(t)\\v(t)\\ \dot u(t)\\\dot v(t)
 \end{smallmatrix}\right)=P\left(\begin{smallmatrix}w_1(t)\\w_2(t)\\ \dot w_1(t)\\\dot w_2(t)
 \end{smallmatrix}\right),\quad t\geq 0,
\ee 
where 
$$
P=\left(\begin{matrix}                            I & 0 & 0 & 0 \\

                                                      -A_2^{-\half}C^* & 0 & 0 &A_2^{-\half} \\
                                                       0 & 0 & I & 0 \\
                                                      0 &  -A_2^{\half}&0 & 0
                                                    \end{matrix}\right).
 $$
Together with \eqref{damped002}, derivation of the equation \eqref{change} leads to the 
coupled systems \eqref{damped u}-\eqref{damped v}. 
The initial data \eqref{outputuv} follows from \eqref{change}.

By Remark \ref{remarque},   $P$ is a bounded and invertible operator from ${\mathcal  H}:=  H_{1,\half}\times H_{2,\half}\times H_1\times H_2 $  to $\mathcal{H} $ with inverse
 $$P^{-1}=\left(
            \begin{matrix}
              I & 0 & 0 & 0 \\
              0 & 0 & 0 & -A_2^{-\half} \\
              0 &0& I & 0 \\

              C^* & A_2^{\half}&0 & 0 \\
            \end{matrix}
          \right).$$
Using $P^{-1}$, the converse in Theorem \ref{transform} can be similarly proved.
 \endproof

The equivalence of the well-posedness of the systems \eqref{damped001}-
\eqref{damped002} and \eqref{damped u}-\eqref{outputuv} can be proved also by using 
their corresponding Cauchy problems. Roughly speaking,  setting $X:=\left(  \begin{smallmatrix}
                w_1 \\
                w_2 \\
                \dot w_1 \\
                \dot w_2 \\
              \end{smallmatrix}
            \right),$ the system \eqref{damped001}-\eqref{damped002}  can be transformed in ${\mathcal  H}$ to the following first order system
 \be\label{cauchy1}\dot X=\mathcal{A}_1 X,\quad X(0)=\left(  \begin{smallmatrix}
                w_1^0 \\
                  w_2^0 \\
                  w_1^1 \\
                w_2^1 \\
              \end{smallmatrix}
            \right),\ee

 where

 \be
  \label{op1b} {\mathcal  A}_1 : {\mathcal
D}({\mathcal  A}_1) \subset {\mathcal  H} \longrightarrow {\mathcal
H}, \, {\mathcal A}_1\left(
\begin{smallmatrix}u_1\\v_1\\u_2\\v_2\end{smallmatrix}\right) = \left(
\begin{smallmatrix}
u_2 \\
v_2\\
- A_1u_1 - \, BB^*u_2  -C v_2\\
-A_2v_1+C^*u_2
\end{smallmatrix}
\right), \ee and $$ {\mathcal  D}({\mathcal  A}_1) :=
 \Bigr\{ (u_1,v_1,u_2,v_2) \in  H_{1,\half} \times H_{2,1}\times H_{1,\half}
 \times H_{2,\half} ,\,\, \, A_1 u_1 +
BB^*u_2  \in H_1\Bigr\}.
$$
 The system \eqref{damped u}-\eqref{outputuv} can be written as
 \be\label{cauchy2}\dot Y=\mathcal{A}_2Y ,\quad Y(0)=\left(  \begin{smallmatrix}
                u^0 \\
                 v^0 \\
                u^1 \\
                 v^1 \\
              \end{smallmatrix}
            \right),\ee
  where

 \be
  \label{op1} {\mathcal  A}_2 : {\mathcal
D}({\mathcal  A}_2) \subset {\mathcal  H} \longrightarrow {\mathcal
H}, \, 
{\mathcal A}_2\left(
\begin{smallmatrix}u_1\\v_1\\u_2\\v_2\end{smallmatrix}\right) = \left(
\begin{smallmatrix}
u_2 \\
v_2\\
- A_1u_1-C(C^*u_1+A_2^{\half} v_1) - \, BB^*u_2  \\
-A_2^{\half}(C^*u_1+A_2^\half v_1)
\end{smallmatrix}
\right), \ee and \be \notag {\mathcal  D}({\mathcal  A}_2) :=
 \Bigr\{ (u_1,v_1,u_2,v_2) \in  H_{1,\half} \times H_{2,\half}\times H_{1,\half}
 \times H_{2,\half} ,\,\, \, A_1 u_1 +
BB^*u_2\in H_1,\,\,C^*u_1+A_2^\half v_1\in H_{2,\half}  \Bigr\}.\ee
For every $\left(  \begin{smallmatrix}
                u_1 \\
                 v_1 \\
                u_2\\
                 v_2 \\
              \end{smallmatrix}
            \right)\in \cD({\mathcal A}_1)$,  we have 
\begin{align*}P\left(  \begin{smallmatrix}
                u_1 \\
                 v_1 \\
                u_2\\
                 v_2 \\
              \end{smallmatrix}
            \right)&=\left(  \begin{smallmatrix}
                u_1 \\
                 -A_2^{-\half}Cu_1 +A_2^{-\half}v_2\\
                u_2\\
                -A_2^\half v_1\\
              \end{smallmatrix}
            \right)
            \\&\\
            &=\left(  \begin{smallmatrix}
                u_1 \\
                 -A_2^{-1}A_2^{\half}Cu_1 +A_2^{-\half}v_2\\
                u_2\\
                -A_2^\half v_1\\
              \end{smallmatrix}
            \right).\end{align*}
Since $C^* u_1+A_2^\half(-A_2^{-\half}Cu_1 +A_2^{-\half}v_2)= v_2 \in H_{2,\half}$ and $A_1 u_1 +
BB^*u_2\in H_1$, we have  $P\cD({\mathcal A}_1)\subset \cD({\mathcal A}_2).$
Using \eqref{change}, we can see that $\mathcal{A}_1=P^{-1} \mathcal{A}_2 P.$

To study the well-posedness and exponential stability of both coupled systems, we write the system \eqref{damped u}-\eqref{outputuv},
 in the product space  $H:=H_1\times H_2$, as  the following second order system
  \BEQ{second}
  \ddot W (t)+A W(t)+\mathcal{B}_0\mathcal{B}_0^*\dot W(t)=0,\quad t\geq 0,
  \EEQ
  \BEQ{second0}W(0)=W^0,\quad \dot W(0)=W^1,\EEQ where 
\be\label{a,b}
 A : {\mathcal
D}( A) \subset {  H} \longrightarrow {
H}, \,  A\left(
\begin{smallmatrix}u\\v\end{smallmatrix}\right) = \left(
\begin{smallmatrix}
A_1u+C(C^*u+A_2^{\half} v)\\
A_2^{\half}(C^*u+A_2^\half v)
\end{smallmatrix}
\right),   \quad \mathcal{B}_0=\left(
                  \begin{smallmatrix}
                    B \\
                    0  \\
                  \end{smallmatrix}
                \right),
 \ee 
with ${\mathcal
D}({  A})=\{(u,v)\in H_{1,1}\times H_{2,\half},\,\,C^*u+A_2^{\half} v\in H_{2,\half}\}$.

 To obtain the well-posedness result,  we need the following lemma which will be also crucial for the rest of this paper.
  \begin{lemma}
  The following assertions  hold.
  \\
  (i)\,\, The operator $A$ is positive self adjoint.\\
    (ii) \,\,$\mathcal{B}_0^* : \cD(A^\half)\longrightarrow U$ is a bounded operator.\\
(iii)\,\,$\mathcal{B}_0 : U\longrightarrow \cD(A^\half)^*:=H_{\frac12}$ is a bounded operator.
  \end{lemma}
  \proof
(i) Let $\left(
              \begin{smallmatrix}
                x \\
                y \\
              \end{smallmatrix}
            \right)\in H_{1,1}\times H_{2,1}.$
            We have

\be\notag <A\left(
              \begin{smallmatrix}
                x \\
                y \\
              \end{smallmatrix}
            \right),\left(
              \begin{smallmatrix}
                x \\
                y \\
              \end{smallmatrix}
            \right)>=  \|A_1^{\half}x\|^2+\langle C^*x+A_2^\half y, C^*x\rangle +\langle C^*x+A_2^\half y, A_2^\half y\rangle.\ee
        Then    \be \label{positive}<A\left(
              \begin{smallmatrix}
                x \\
                y \\
              \end{smallmatrix}
            \right),\left(
              \begin{smallmatrix}
                x \\
                y \\
              \end{smallmatrix}
            \right)>=\|A_1^{\half}x\|^2+\|A_2^{\half}y+C^*x\|^2>0.\notag\ee

Thus, $A$ is a symmetric positive  operator.
For every  $(f,g)\in H$,  the  solution $(u,v)\in \cD(A)$ of the system
\be\label{eqp1}
A_1u+C(C^*u+A_2^{\half} v)=f,\notag\ee
\be\label{eqp2}
A_2^{\half}(C^*u+A_2^\half v)=g\notag\ee
is given by
\be u=A_1^{-1}(f-CA_2^{-\half}g),\quad v=A_2^{-1}g-A_2^{-\half}C^*A_1^{-1}(f-CA_2^{-\half}g).\notag\ee
It is clear that $(u,v)\in H_{1,\half}\times H_{2,\half}.$ Since $C^*u+A_2^{\half} v=A_2^{-\half}g\in H_{2,\half},$ we have $(u,v)\in \cD(A).$ 
Thus, the operator $A$ is invertible.
Consequently, $A$ is a  positive self adjoint operator.

(ii) Let $\left(
              \begin{smallmatrix}
                x \\
                y \\
              \end{smallmatrix}
            \right)\in \cD(A^\half)$. We have $\mathcal{B}_0^*\left(
              \begin{smallmatrix}
                x \\
                y \\
              \end{smallmatrix}
            \right)= B^*x.$ Since  $B^*$ is a bounded operator from $H_{1,\half}$ to $U$, there exists a constant $c>0$ such that $\|B^*x\|_U\leq c\|A_1^\half x\|_{H_1}$. Thus,
             $$
             \|\mathcal{B}_0^*\left(
              \begin{smallmatrix}
                x \\
                y \\
              \end{smallmatrix}
            \right)\|_U\leq c [\|A_1^{\half}x\|^2+\|A_2^{\half}y+C^*x\|^2],
$$
and thus  the operator $\mathcal{B}_0^* : \cD(A^\half)\longrightarrow U$ is  bounded. 
The assertion (iii) follows from (ii).
 \endproof
  As a consequence of the above lemma  we have the following well-posedness result.
\begin{proposition}\label{propcond}
Assume that \eqref{as00}, \eqref{as0} and \eqref{as1} hold.
  Then,
the system \be \label{wave}
  \ddot W (t)+A W(t)+\mathcal{B}_0\mathcal{B}_0^*\dot W(t)=0,\quad t\geq0,
  \ee
  \be \label{initialesbb} W(0)=W^0,\,\dot W(0)=W^1\ee is well-posed in the energy space $\cD(A^{\half})\times H.$
\end{proposition}

    Using Theorem \ref{transform},  Proposition \ref{propcond} and the regularity results in \cite{at}, we have the following  well-posedness and regularity result of the coupled systems  \eqref{damped001}-\eqref{output00} .
\begin{proposition}\label{propcond1}
Assume that \eqref{as00}, \eqref{as0} and \eqref{as1} hold.
  Then, the  system \eqref{damped001}-\eqref{output00} is well-posed, i.e., 
\begin{itemize}
                                                               \item[(i)] for $(w_1^0,w_2^0,w_1^1,w_2^1)\in \cD(\cA_1)$, the problem \eqref{damped001}-\eqref{output00} admits a unique
solution
$w_i \in C^1([0,T];H_{i,\half}) \cap C^2([0,T];H_i), i = 1, 2,$
                                                               \item[(ii)] for $(w_1^0,w_2^0,w_1^1,w_2^1)\in \cH$,
$w_i \in C([0,T];H_{i,\half}) \cap C^1([0,T];H_i), i = 1, 2,$ and $B^*w_1(\cdot) \in H^1(0, T;U).$

                                                             \end{itemize}
\end{proposition}
\begin{remark}
The well-posedness of \eqref{damped001}-\eqref{output00} can be also obtained directly by proving that the operator 
$\cA_1$ satisfies the conditions of Lumer-Phillips theorem, see \cite{engel}.
\end{remark}
\subsection{Transfer function}
To characterize the stabilization of system \eqref{damped001}-\eqref{output00} we need  the following lemma.
\begin{lemma}
Assume that \eqref{as00}, \eqref{as0} and \eqref{as1} hold.
  Then, the following results hold.
\begin{itemize}
  \item[(i)] The operator $[\lambda^2+A_1+\lambda^2C(\lambda^2+A_2)^{-1}C^*]$ is invertible from $H_{1,\half}$ to $H_{1,-\half}$.
  \item[(ii)] The function defined  by $$G(\lambda)=\lambda B^*[\lambda^2+A_1+
\lambda^2C(\lambda^2+A_2)^{-1}C^*]^{-1}B,\quad \lambda >0,$$ is the transfer function of
both systems \eqref{damped001}-\eqref{damped002} and \eqref{wave}-\eqref{initialesbb}.
\end{itemize}
\end{lemma}
\proof
\begin{itemize}
  \item[(i)] Let $y\in H_{1,-\half}.$ 
Consider in $H_{1,\half}$ the equation \be\label{lax}[\lambda^2+A_1+\lambda^2
C(\lambda^2+A_2)^{-1}C^*]x=y.\ee  For  every $\zeta\in H_{1,\half}$,  we have
$$\left<[\lambda^2+A_1+\lambda^2C(\lambda^2+A_2)^{-1}C^*]x,\zeta\right>=\left<y,\zeta\right>$$
which can be written as
$$
\lambda^2\left<x,\zeta\right>+\left<A_1^\half x,A_1^\half\zeta\right>+
\left<\lambda(\lambda^2+A_2)^{-\half}C^*x,\lambda(\lambda^2+A_2)^{-\half}C^*\zeta\right>=
\left<y,\zeta\right>=:\Lambda(x,\zeta).
$$
 Since  $\Lambda$ is a
bilinear coercive form on $H_{1,\half}$, the Lax-Milgram theorem leads to
the existence and uniqueness of the solution $x$ to  the equation
\eqref{lax}, and thus the claim follows.
  \item[(ii)] We compute first the transfer function  of \eqref{wave}-\eqref{initialesbb}.
Setting $Z: =\left(\begin{smallmatrix} W\\
\dot W
\end{smallmatrix}\right)$, the open loop system associated to  \eqref{wave}-\eqref{initialesbb} can be transformed  to the following controlled first order system in the  energy space $\cD(A^\half)\times H$
\begin{align}\label{csecond}
\dot Z (t)&= \mathcal{A}_2^0 Z(t) + \mathcal{B}u(t), \quad t\geq 0,
\\
Z(0)&=\left(\begin{smallmatrix} W^0\\
 W^1
\end{smallmatrix}\right),
\label{initisecond}
\end{align}
with $\mathcal{A}_2^0=\left(
           \begin{matrix}
             0 & I \\
             -A & 0 \\
           \end{matrix}
         \right),$ $\cD(\cA_2^0)=\cD(A)\times\cD(A^\half),$\quad and  $\mathcal{B}=\left(
              \begin{matrix}
                \mathcal{B}_0 \\
                0\\
              \end{matrix}
            \right).$

Let $\left(
              \begin{smallmatrix}
                f \\
                g \\
              \end{smallmatrix}
            \right)\in \cD(A^\half)\times H.$  We look  for $\left(
              \begin{smallmatrix}
                x \\
                y \\
              \end{smallmatrix}
            \right)\in \cD(A)\times \cD(A^\half)$ such that
 \be \label{2resolvent}(\lambda  -\mathcal{A}_2^0) \left(
              \begin{smallmatrix}
                x \\
                y \\
              \end{smallmatrix}
            \right)=\left(
              \begin{smallmatrix}
                f \\
                g \\
              \end{smallmatrix}
            \right).\ee We have

$\begin{array}{ccl}
\eqref{2resolvent}&\Longleftrightarrow &\begin{cases}\lambda x - y=f\\
\lambda y+Ax=g
\end{cases}\\
      &\Longleftrightarrow &\begin{cases}x=\lambda (\lambda ^2+A)^{-1}f +(\lambda ^2+A)^{-1}g\\
y=(\lambda^2 (\lambda ^2+A)^{-1}-I)f+\lambda (\lambda ^2+A)^{-1}g,
\end{cases}
\end{array}
$

and thus $$(\lambda  -\mathcal{A}_2^0)^{-1}=\left(
                                      \begin{smallmatrix}
                                        \lambda (\lambda ^2+A)^{-1} & (\lambda ^2+A)^{-1} \\
                                        (\lambda^2 (\lambda ^2+A)^{-1}-I) & \lambda (\lambda ^2+A)^{-1} \\
                                      \end{smallmatrix}
                                    \right).
$$
 The transfer function $G_2(\lambda):= \mathcal{B}^*(\lambda-\mathcal{A}_2^0)^{-1}\mathcal{B}$ of the system 
\eqref{csecond}-\eqref{initisecond} is then
\be \notag G_2(\lambda)=\left(\begin{matrix}\mathcal{B}_0^* &0\end{matrix}\right)\left(
              \begin{smallmatrix}
                \lambda (\lambda ^2+A)^{-1} \mathcal{B}_0 \\
                (\lambda^2 (\lambda ^2+A)^{-1}-I)\mathcal{B}_0 \\
              \end{smallmatrix}
            \right)=\lambda \mathcal{B}_0^*
                 (\lambda ^2+A)^{-1} \mathcal{B}_0 .
\ee

Easy computation  leads to $$(\lambda^2+A)^{-1}=\left(
                                                            \begin{matrix}
                                                              \Gamma && -\Gamma C A_2^\half
                                                              (\lambda^2+A_2)^{-1} \\
                                                              -(\lambda^2+A_2)^{-1}A_2^\half
                                                              C^*\Gamma &&(\lambda^2+A_2)^{-1}[I+A_2^\half C^*\Gamma C A_2^\half(\lambda^2+A_2)^{-1} ] \\
                                                            \end{matrix}
                                                         \right),$$
where $\Gamma:=[\lambda^2+A_1+\lambda^2C(\lambda^2+A_2)^{-1}C^*]^{-1} .$
Consequently,   \be\label{2transfer} G_2(\lambda)=
\lambda B^*[\lambda^2+A_1+\lambda^2C(\lambda^2+A_2)^{-1}C^*]^{-1}B ,\quad \forall \lambda >0.\notag\ee

Let ${\mathcal  A}_1^0 : {\mathcal
D}({\mathcal  A}_1^0) \subset {\mathcal  H} \longrightarrow {\mathcal
H}, \, {\mathcal A}_1^0\left(
\begin{array}{ccc}u_1\\v_1\\u_2\\v_2\end{array}\right) = \left(
\begin{array}{l}
u_2 \\
v_2\\
- A_1u_1   -C v_2\\
-A_2v_1+C^*u_2
\end{array}
\right)$   be the generator of the open loop system associated to \eqref{damped001}-\eqref{output00}. 
Since $\mathcal{A}_1^0=P^{-1} \mathcal{A}_2^0 P$, we have
\be\label{resolv}
(\lambda -\mathcal{A}_1^0)^{-1}= P^{-1}(\lambda -\mathcal{A}_2^0)^{-1}P,\quad \forall \lambda >0.\notag\ee
 Since $\mathcal{B}_0^*P^{-1}=\mathcal{B}_0^*$ and $P\mathcal{B}_0=\mathcal{B}_0$,  we have 
$$G_1(\lambda):=\mathcal{B}_0^*(\lambda -\mathcal{A}_1^0)^{-1}\mathcal{B}_0=
\mathcal{B}_0^*(\lambda -\mathcal{A}_2^0)^{-1}\mathcal{B}_0=G_2(\lambda).$$

\end{itemize}
\subsection{Stabilization }
In order to characterize the stabilization of the coupled systems without delay, we give some energy equivalences.
\begin{proposition}
Assume that \eqref{as00}, \eqref{as0}, \eqref{as1}
 and \eqref{condition} hold. Then,
 $$\begin{array}{ccc}&\mathcal{E}(t):=\dfrac{1}{2}\left( \|A^{\half}     \left(  \begin{smallmatrix}
                u \\
                v \\
              \end{smallmatrix}
            \right)\|^2+\|     \left(  \begin{smallmatrix}
                \dot u \\
                \dot v \\
              \end{smallmatrix}
            \right)\|^2_{H_1\times H_2}\right)&\\
            &\asymp&\\
  &\widetilde{E}(t):=\dfrac{1}{2}\left(\|A_1^{\half}u\|^2+\|A_2^{\half}v\|^2+\|C^*u\|^2+\|     \left(  \begin{smallmatrix}
                \dot u \\
                \dot v \\
              \end{smallmatrix}
            \right)\|^2_{H_1\times H_2}\right)&\\
            &\asymp&\\
            &E(t):= \dfrac{1}{2}\left(\|A_1^{\half}w_1\|^2+\|A_2^{\half}w_2\|^2+\|     \left(  \begin{smallmatrix}
                \dot w_1 \\
                \dot w_2\\
              \end{smallmatrix}
            \right)\|^2_{H_1\times H_2}\right),&
            \end{array}$$
for every solutions $(w_1,w_2)$ and $(u,v)$ of \eqref{damped001}-\eqref{output00} and 
\eqref{damped u}-\eqref{outputuv}, respectively.
\end{proposition}
From this follows immediately the following corollary.

 \begin{corollary}\label{cora}
Assume that \eqref{as00}, \eqref{as0}, \eqref{as1}
 and \eqref{condition} hold. Then,
The exponential stabilities of the three  systems \eqref{damped001}-\eqref{output00}, 
\eqref{damped u}-\eqref{outputuv}, and  \eqref{wave}-\eqref{initialesbb} are equivalent.
\end{corollary}
Using the characterization of stabilization of second order equation in \cite{at}, we have the following result.
\begin{theorem}  \label{theor} Assume that \eqref{as00}, \eqref{as0}, \eqref{as1}
 and \eqref{condition} hold  and for  a fixed $\gamma>0$
\begin{equation}\label{BOUNDEDTRANSFERb}
\sup_{{\rm Re}\lambda =\gamma}\left\|\lambda B^{*}\left[\lambda^2I+A_1+\lambda^2C(\lambda^2+A_2)^{-1}C^*\right]^{-1}
B \right\|_{{\mathcal  L}(U)}<\infty \m.
\end{equation}
The system \eqref{wave}-\eqref{initialesbb} is exponentially stable in $\cD(A^{\half})\times H$ if and only if there exists a constant $c>0$ such that \begin{align}\label{io1} c\int_0^T\|(B^*\phi_1)'\|_U^2dt\geq\|A^{\half}\left(  \begin{smallmatrix}
                 \phi_1(0) \\
                 \psi_1(0) \\
              \end{smallmatrix}
            \right)\|^2+\|     \left(  \begin{smallmatrix}
                   \dot \phi_1(0) \\
                \dot \psi_1(0) \\
               \end{smallmatrix}
            \right)\|^2_{H_1\times H_2}, \end{align} where $(\phi_1,\psi_1)$ is a solution of the following system
\be\label{conserv1}\ddot \phi_1 + (A_1+CC^*)\phi_1 +CA_2^{\half}\psi_1=0\ee
\be\label{conserv2}\ddot \psi_1 +A_2 \psi_1+A_2^{\half}C^*\phi_1=0.\ee
\end{theorem}
As a consequence of the above theorem, we have the following result.
\begin{theorem}\label{main1} Assume that \eqref{as00}, \eqref{as0}, \eqref{as1},
 \eqref{condition}, and \eqref{BOUNDEDTRANSFERb} hold.  Then the following assertions are equivalent.
\\
  (i)\, The system \eqref{damped001}-\eqref{output00} is exponentially stable in $\mathcal{H}.$
\\
  (ii)\, There exists a constant $c>0$ such that \begin{align} 
\label{io2b} c\int_0^T\|(B^*\phi)'\|_U^2dt\geq\|A_1^{\half}\phi(0)\|^2+\|A_2^{\half} \psi(0)\|^2+\|C^*\phi(0)\|^2+\|     \left(  \begin{smallmatrix}
                \dot \phi(0) \\
                \dot \psi(0) \\
              \end{smallmatrix}
            \right)\|^2_{H_1\times H_2}.
\end{align}
  (iii)\,  There exists a constant $c>0$ such that \begin{align} 
\label{io2} c\int_0^T\|(B^*\phi)'\|_U^2dt\geq\|A_1^{\half}\phi(0)\|^2+\|A_2^{\half} \psi(0)\|^2+\|     \left(  \begin{smallmatrix}
                \dot \phi(0) \\
                \dot \psi(0) \\
              \end{smallmatrix}
            \right)\|^2_{H_1\times H_2},\end{align}
 where $(\phi,\psi)$ is a solution of the following conservative adjoint system
\be\label{conserv11}\ddot \phi + A_1\phi +C\dot\psi=0,\ee
\be\label{conserv21}\ddot \psi +A_2 \psi-C^*\dot\phi=0.\ee
\end{theorem}
\proof
From Corollary \ref{cora} and Theorem \ref{theor}, the assertion (i) is equivalent to the observability inequality \eqref{io1}. 
To show \eqref{io2b} and \eqref{io2}, let $(\phi_1,\psi_1)$ be a solution of \eqref{conserv1}-\eqref{conserv2}. Then 
 $\phi=\phi_1$ and $\psi=A_2^{-\half}\dot \psi_1$ satisfy  \eqref{conserv11}-\eqref{conserv21}.
The observability inequality \eqref{io1} becomes
 $$c\int_0^T\|(B^*\phi)'\|_U^2dt\geq\|A_1^{\half}\phi(0)\|^2+\|A_2^{\half} \psi(0)\|^2+\|C^*\phi(0)\|^2+\|     \left(  \begin{smallmatrix}
                \dot \phi(0) \\
                \dot \psi(0) \\
              \end{smallmatrix}
            \right)\|^2_{H_1\times H_2}+2\Re e<\phi(0),C\dot\psi(0)>.$$

Since by \eqref{condition}, 
 $$|<\phi,C\dot\psi(0)>|<\delta\left(\|A_1^{\half}\phi(0)\|^2+\|\dot\psi(0)\|^2+\|C^*\phi(0)\|^2\right)$$
with $\delta\in[0,\half),$
 the inequality \eqref{io1} can be written as
$$c\int_0^T\|(B^*\phi)'\|_U^2dt\geq\|A_1^{\half}\phi(0)\|^2+\|A_2^{\half} \psi(0)\|^2+\|C^*\phi(0)\|^2+\|     \left(  \begin{smallmatrix}
                \dot \phi(0) \\
                \dot \psi(0) \\
              \end{smallmatrix}
            \right)\|^2_{H_1\times H_2}, 
$$
which is exactely the inequality \eqref{io2b}. Now from the assumption $H_{1,\half}\hookrightarrow \cD(C^*)$ 
follows the inequality \eqref{io2}. The converse can be  shown in the same way.
\endproof
\section{Coupled second order systems with delay}
Consider in this section the following coupled systems with delay \BEQ{damped12} \ddot
w_1(t) + A_1w_1(t) +  \alpha_1 \, BB^*\dot w_1(t)+ \alpha_2 \,
BB^*\dot w_1(t-\tau)+ C \dot w_2(t) = 0, \quad t\geq 0, \EEQ
\BEQ{damped22} \ddot w_2(t) + A_2w_2(t) -C^* \dot
w_1(t) = 0, \quad t\geq 0,\EEQ
 \BEQ{output2} w_i(0) \m=\m w_i^0,  \m
\dot w_i(0)=w_i^1,\,i=1,2, \,\, \dot w_1(s) = f_0(s), \,\,s\in
(-\tau,0), \EEQ where $\tau
> 0$ is the time delay, $\alpha_1$ and $\alpha_2$ are positive real
numbers, and the initial data $(w_1^0,w_1^1,w_2^0,w_2^1,f_0)$
belongs to a suitable space.

 Using the same method as in the coupled systems without delay, the system \eqref{damped12}-\eqref{damped22} can be transformed to the following one
 \BEQ{damped ud}
 \ddot u(t) +(A_1+CC^*)u(t)+ \alpha_1BB^*\dot u(t)+\alpha_2BB^*\dot u(t-\tau)+ CA_2^{\half}v(t)=0 , \quad t\geq 0,
 \EEQ
 \BEQ{damped vd}
 \ddot v(t)+ A_2 v(t) + A_2^{\half}C^*u(t)=0 , \quad t\geq 0.
 \EEQ
\BEQ{outputuvd}
 u(0)=u^0, \dot u(0)=u^1, v(0)=v^0, \dot v(0)=v^1
 \,\, \dot u(s) = f_0(s), \,\,s\in
(-\tau,0),\EEQ
with
$
 u^0:=w_1^0,  u^1:=w_1^1, v^0:=A_2^{-\half} w_2^1-A_2^{-\half}C^*w_1^0,  v^1:=-A_2^{\half}w_2^0.
 $
This   system can be written in the space $H=H_1\times H_2$ under the following second order system with 
delay
  \BEQ{secondd}
  \ddot W (t)+A W(t)+\alpha_1\mathcal{B}_0\mathcal{B}_0^*\dot W(t)+\alpha_2\mathcal{B}_0\mathcal{B}_0^*\dot W(t-\tau)=0,
  \EEQ
  \BEQ{initial data}
  W(0)=W^0,\,\,\dot W(0)=W^1,\,\,\dot W(s)= \left(
              \begin{smallmatrix}
                f_0(s)\\
                0\\
              \end{smallmatrix}
            \right),\,\,s\in
(-\tau,0),
  \EEQ with $A$ and $ \mathcal{B}_0$ are defined in the previous section.
Let $E_{1,\half}$
be the topological supplement of $ker B^*$ in $H_{1,\half}$ and $P_2$ its associated projection. It is clear that $E_{1,\half}\times \{0\}$ is the topological supplement of $ker \cB_0^*$ in $H_{\half}=H_{1,\half}\times H_{2,\half}$ and the associated projection $\mathcal{P}_2$  is given by
$\mathcal{P}_2W^0=\left(
                  \begin{matrix}
                    P_2u^0 \\
                    0  \\
                  \end{matrix}
                \right).
                $
As in \cite{abmb}, the second order equation with delay \eqref{secondd}-\eqref{initial data} can
 be written  as the  Cauchy problem 
\BEQ{cauchy}
  \dot {\widetilde{Z}} (t)=\cA_d \widetilde{Z}(t) ,\,\,\, t\geq0, \quad \widetilde{Z}(0)=\left(\begin{smallmatrix}W^0\\W^1\\P_2 f_0 
\end{smallmatrix}\right)
\ee 
in the Hilbert space $ H_\half\times H\times L^2(-\tau,0;\mathcal{P}_2H_\half)$ which  can be identified with 
$\widetilde{\cH}:= H_\half\times H\times L^2(-\tau,0; E_{1,\half})$, 
where $\widetilde{Z}= \left (\begin{smallmatrix}u\\v\\ \dot u\\ 
\dot v\\z\end{smallmatrix}\right)$, $z(t,\theta)=P_2\dot u(t+\theta),\,\,\, \theta \in (-\tau,0)$ and
\be
  \label{op} {\mathcal A}_d\left(
\begin{smallmatrix}u_1\\v_1\\u_2\\v_2\\z
\end{smallmatrix}\right) = \left(
\begin{smallmatrix}
u_2 \\
v_2\\
- A_1u_1-C(C^*u_1+A_2^{\half} v_1) - \alpha_1 \, BB^*u_2 - \alpha_2 \, BB^*z(-\tau) \\
-A_2^{\half}(C^*u_1+A_2^{\half} v_1)\\
 \partial_\theta z
\end{smallmatrix}
\right), \ee 
 with \be\notag {\mathcal  D}({\mathcal  A}_d) :=
 \label{d}\notag\Bigr\{ (u_1,v_1,u_2,v_2,z) \in  H_{1,\half}\times H_{2,\half} \times H_{1,\half}
 \times H_{2,\half} \times
H^1(-\tau,0;E_{1,\half}),\ee \be \notag \, A_1 u_1 + \alpha_1 \,
BB^*u_2 + \alpha_2 BB^*z(-\tau) \in H_1,\;C^*u_1+A_2^{\half} v_1\in H_{2,\half},\; z(0)=P_2u_2 \Bigr\}.\ee
Assuming $\alpha_2 \leq \alpha_1,$ we introduce in $\widetilde{\cH}$ the new inner product
$$
\left<\left(
\begin{array}{ccc}U_1\\U_2\\z_1\end{array}\right),\left(
\begin{array}{ccc}V_1\\V_2\\z_2\end{array}\right)\right>=
\left<U_1,V_1\right>_{H_{\half}}+\left<U_2,V_2\right>_{H}+
\frac{\xi}{\tau}\int_{-\tau}^0\left<B^*z_1(\theta),B^*z_2(\theta)\right>_U
d\theta,
$$
where  $\xi$ is a constant satisfying
\begin{equation}\label{alfa}
\tau \alpha_2 \leq   \xi \leq \tau (2 \alpha_1 - \alpha_2).
\end{equation}
It can be seen easily that $\widetilde{\cH}$ endowed with this inner
product is a Hilbert space, and   its associated norm  is equivalent
to the canonical norm of $\widetilde{{\cal H}}$. Now,  we are in the position to use the results in 
\cite{abmb} to \eqref{damped12}-\eqref{output2}, and deduce first its well-posedness.
To characterize the  stabilization, we introduce
the following delay energy functions
$$
E_d(t) :=
$$
\be
\label{energyb}
\half \, \left\|(w_1(t), w_2(t))\right\|^2_{H_{1,\half}\times H_{2,\half}}
+\half \left\|(\dot w_1(t),\dot w_2(t))\right\|^2_{H_1\times H_2} +
\frac{\xi}{2} \, \int_{-\tau}^0
\left\|B^* \dot w_1(t+\theta)\right\|^2_U \, d \theta, \quad t\geq 0,\notag
\ee
and
\be
\label{energy}
 \widetilde{E}_{d}(t) :=
\half \, \left\|(u(t), v(t))\right\|^2_{H_{\half} }
+\half \left\|(\dot u(t),\dot v(t))\right\|^2_{{H} } +
\frac{\xi}{2} \, \int_{-\tau}^0 \left\|B^* \dot
u(t+\theta)\right\|^2_U \, d \theta, \quad t\geq 0. 
\notag
\ee
Under the assumption \eqref{condition}, $E_d(t)$ and $\widetilde{E}_{d}(t)$ are equivalent.

By our result \cite[Theorem 1.1]{abmb},  Theorem \ref{main1} yields the following main result.
\begin{theorem} \label{obsdel}
Assume that   \eqref{as00}, \eqref{as0}, \eqref{as1},  \eqref{condition} and   \eqref{BOUNDEDTRANSFERb} hold and that  
$\alpha_2<\alpha_1$.
           
Then the following  assertions are equivalent.
\begin{enumerate}
  \item There are constants $\omega,C>0$ such that the system
\rfb{damped12}-\rfb{output2} satisfies the exponential decay
$$ E_d(t)\leq Ce^{-\omega t}E_d(0)\,\,  \mbox{for all}\,\,
(w_1^0,w_2^0,w_1^1,w_2^1,f_0)\in \mathcal{D}(\cA_d).$$

%$$\mathcal{D}=\{(u_0,v_0,u_1,v_1,z)\in H_{1,\half}\times H_{2,1}\times H_{1,\half}\times  H_{2,\half}\times H^1(-\tau,0,E_{1,\half}),$$$$ A_1u_0+\alpha_1 BB^* u_1+\alpha_2 BB^*z(-\tau)\in H_1, \,z(0)=u_1\}.
%$$
  \item There exist $T,c > 0$ such that \begin{align} \label{iod2} c\int_0^T\|(B^*\phi)'\|_U^2dt\geq\|A_1^{\half}\phi(0)\|^2+\|A_2^{\half} \psi(0)\|^2+\|     \left(  \begin{smallmatrix}
                \dot \phi(0) \\
                \dot \psi(0) \\
              \end{smallmatrix}
            \right)\|^2_{H_1\times H_2}\notag\end{align}
 for every solution  $(\phi,\psi)$  of the  conservative adjoint  system
\be\notag \ddot \phi + A_1\phi +C\dot\psi=0,\ee
\be\notag \ddot \psi +A_2 \psi-C^*\dot\phi=0.\ee
\end{enumerate}
 \end{theorem} 
\section{Applications}
\subsection{First example : Dirichlet boundary conditions}
Consider the following  coupled wave equations
\begin{align}\label{Ap1b}
&\ddot w_1(t,x)-\frac{\partial^2 w_1}{\partial x^2}(t,x)+\alpha_1 \dot w_1(t,\xi)\delta_\xi+\alpha_2 \dot
 w_1(t-\tau,\xi)\delta_\xi+\beta \frac{\partial \dot w_2}{\partial x}(t,x)=0, (t,x)\in(0,\infty)\times ]0,1[,\\
&\ddot w_2(t,x)-\frac{\partial^2 w_2}{\partial x^2}(t,x)+\beta \frac{\partial \dot w_1}{\partial x}(t,x)=0,\quad (t,x)\in(0,\infty)\times ]0,1[,\\
&w_i(t,0)=w_i(t,1)=0,\quad  t\in (0,\infty),\; i=1,2, \\
&w_i(0,x)=w_i^0(x),\,\,\dot w_i(0,x)=w_i^1(x),\,\,\dot w_1(s,x)=f_0(s,x), -\tau\leq s < 0,x\in ]0,1[, i=1,2,\label{Ap1c}
\end{align}
with $\xi\in (0,1)$, $\beta>0$ and $0<\alpha_2<\alpha_1$. To put this control system into the framework of this paper, consider the spaces  $H_1=H_2=L^2(0,1)$ and the operators
$
A_1=A_2=-\frac{d^2}{dx^2},
$
with the domain   $\cD(A_1)=\cD(A_2)=H^2(0,1)\cap H^1_0(0,1)$
  which are obviously  self-adjoint positive  operators. In this case, the domains of the  fractional power  operators are given by
$$
\cD(A^\frac12_1)=\cD(A^\frac12_2)=H^1_0(0,1).
$$The operator $B$ and its adjoint $B^*$  are given by
$$
Bk= k\delta_\xi,\;\; k\in\mathbb{R}, \;\;B^*\varphi=\varphi(\xi), \;\; \varphi\in H^{1}_0(0,1)
$$ and  finally
$$
C=\beta\frac{d}{dx}, \mbox{ with } \cD(C)=H^1_0(0,1).
$$
It is clear that  $B^*:H^1_0(0,1)\rightarrow \mathbb{R}$ is  bounded and $C^*= - \beta\frac{d}{dx}$  with
\be H^1_0(0,1)\hookrightarrow \cD(C^*)= H^1(0,1).\ee
Now assume that $\beta<1$, then, with a simple integration  by parts,  the condition \eqref{condition} is satisfied  with constant $\delta=\frac{\beta}{2}$ .
Let us now check  the assumption \eqref{BOUNDEDTRANSFERb}.  Since in this example 
 $A_1=A_2$, we can easily see  that
$$
\left[\lambda^2I+A_1+\lambda^2C(\lambda^2+A_2)^{-1}C^*\right]^{-1}=
\frac12\left[\lambda^2I+A_1+\lambda C\right]^{-1}
+\frac12\left[\lambda^2I+A_1+\lambda C^*\right]^{-1}.
$$
Thus, we have the following decomposition of the transfer function
$$
H(\lambda)=\frac\lambda2B^*\left[\lambda^2I+A_1+\lambda C\right]^{-1}B+\frac\lambda2B^*
\left[\lambda^2I+A_1+\lambda C^*\right]^{-1}B:= H_1(\lambda)+H_2(\lambda).
$$
 For every $k\in \R$, the function 
$$
\psi:=\left[\lambda^2I+A_1+\lambda C\right]^{-1}Bk
$$
 satisfies
\begin{align}
     & \lambda^2\psi(x)-\frac{d^2 \psi}{dx^2}(x)+ \lambda\beta\frac{d \psi}{dx}(x)=0,  x\in (0,\xi)\cup(\xi,1) \label{e1}\\
   & \psi(0)=\psi(1)=0, \label{e2} \\
&[\psi]_\xi=0, \left[\frac{d \psi}{dx}\right]_\xi=k \label{e3},
\end{align}
where we denote by $[g]_\xi$ the jump of the function $g$ at the point $\xi$.
The solutions $r_1, r_2$  of the equation $r^2-\beta\lambda r-\lambda^2=0$  are 
$
\frac\lambda2(\beta\pm\sqrt{\beta^2+4})
$.
Hence,  the solution of \eqref{e1}-\eqref{e2} is 
$$
\psi(x)=\left\{
  \begin{array}{ll}
    A\left(e^{r_1x}-e^{r_2x}\right), & x\in(0,\xi) \\
    B\left(e^{r_1(x-1)}-e^{r_2(x-1)}\right), & x\in(\xi,1),
  \end{array}
\right.
$$
and  \eqref{e3} yields 
$$
\psi(x)=\frac{k e^{- \lambda \beta \xi}}{\lambda\sqrt{\beta^2+4}}\left\{
  \begin{array}{ll}
    \frac{e^{r_1(\xi - 1)}-e^{r_2(\xi -1)}}{e^{-r_2}-e^{-r_1}}\left(e^{r_1x}-e^{r_2x}\right), & x\in(0,\xi) \\
    \frac{e^{r_1\xi}-e^{r_2\xi}}{e^{-r_2}-e^{-r_1}}\left(e^{r_1(x-1)}-e^{r_2(x-1)}\right), & x\in(\xi,1).
  \end{array}
\right.
$$
Consequently
$$
H_1(\lambda)=\frac{e^{- \lambda \beta \xi}}{2\sqrt{\beta^2+4}}\frac{e^{r_1(\xi -1)}-e^{r_2(\xi-1)}}{e^{-r_2}-e^{-r_1}}\left(e^{r_1\xi}-e^{r_2\xi}\right)
$$
and then, for every $\gamma>0$, we have
$$
\sup_{Re\lambda=2\gamma}|H_1(\lambda)|\leq\frac{1}{\sqrt{\beta^2+4}}\frac{\cosh(\gamma\sqrt{\beta^2+4}(1-\xi))}{\sinh(\gamma\sqrt{\beta^2+4})}\cosh(\gamma\sqrt{\beta^2+4}\xi).
$$
By similar calculus, we have the boundedness of $H_2$, and thus the assumption \eqref{BOUNDEDTRANSFERb} is satisfied

Now, consider the conservative adjoint system
\begin{align}\label{Ap1}
&\frac{\partial^2 \phi}{\partial t^2}(t,x)-\frac{\partial^2 \phi}{\partial x^2}(t,x)+\beta \frac{\partial^2 \psi}{\partial x \partial t}(t,x)=0,\quad (t,x)\in(0,\infty)\times ]0,1[,\\
&\frac{\partial^2 \psi}{\partial t^2}(t,x)-\frac{\partial^2 \psi}{\partial x^2}(t,x)+\beta \frac{\partial^2 \phi}{\partial x \partial t}(t,x)=0,\quad (t,x)\in(0,\infty)\times ]0,1[,\\
&\phi(t,0)=\psi(t,0)=\phi(t,1)=\psi(t,1)=0,\quad  t\in (0,\infty), \\
&\phi(0,x)=\phi^0(x),\,\,\dot \phi(0,x)=\phi^1(x),\,\, \psi(0,x)=\psi^0(x),\,\,\dot \psi(0,x)=\psi^1(x),\,\, x\in ]0,1[.
\end{align}

Consider the initial conditions as follows
\begin{align*}
\phi^0(x)=\sum_{n \in \mathbb{Z}^*} a_n \, \cos\left(\frac{n\beta \pi}{\sqrt{\beta^2 + 4}} x \right) \, \sin(n \pi x), \phi^1(x)=\sum_{n \in \mathbb{Z}^*} \lambda_n \, a_n \, \cos\left( \frac{n\beta \pi}{\sqrt{\beta^2 + 4}}x \right) \, \sin(n \pi x)\\
\psi^0(x)=\sum_{n \in \mathbb{Z}^*} a_n \, \sin\left(\frac{n\beta \pi}{\sqrt{\beta^2 + 4}}x \right) \, \sin(n \pi x), \psi^1(x)=\sum_{n \in \mathbb{Z}^*} \lambda_n \,  a_n \, \sin\left( \frac{n \beta \pi}{\sqrt{\beta^2 + 4}}x \right)) \, \sin(n \pi x)
\end{align*}
with $(\lambda_n \, a_n)$ are in $l^2(\mathbb{C}),$ where $\lambda_n = \frac{2i n \pi}{\sqrt{\beta^2 + 4}}, \, \forall \, n \in \mathbb{Z}^*.$

By standard technics, we obtain
$$
\phi(t,x) =\sum_{n \in \mathbb{Z}^*}
a_n \, e^{\lambda_n t} \, \cos\left(\frac{n\beta \pi}{\sqrt{\beta^2 + 4}}x \right) \,
\sin(n \pi x)
$$
and then,
\begin{align*}
\frac{\partial \phi}{\partial t}(t,\xi)=  \sum_{n \in \mathbb{Z}^*} \lambda_n \, a_n \, e^{\lambda_n t}\, \cos\left( \frac{n \beta \pi}{\sqrt{\beta^2 + 4}} \xi \right) \, \sin(n \pi \xi).
\end{align*}
Now, by the Ingham's inequality, for any  $T> \sqrt{\beta^2+4}$ we have
\be
\label{obsineq1}
\int_0^T\left|\frac{\partial \phi}{\partial t}(t,\xi)\right|^2 dt
\asymp \sum_{n \in \mathbb{Z}^*} |\lambda_n|^2 \, |a_n|^2 \, 
\left|\cos\left( \frac{n\beta \pi}{\sqrt{\beta^2 + 4}}\xi \right) \right|^2\, \left| \sin (n \pi \xi) \right|^2,
\ee
%we can prove as above that for $T > \sqrt{\beta^2 + 4}$ there exists a constants $C = C_{T,\xi,\beta} > 0$ such that
%\be
%\label{obsineq2}
%\int_0^T\left|\frac{\partial \phi}{\partial t}(t,\xi)\right|^2 dt
%\geq C_{T,\xi,\beta} \sum_{n \in \mathbb{Z}^*} |\lambda_n|^2 \, |a_n|^2 \, \left|\cos\left( \frac{n\beta %\pi}{\sqrt{\beta^2 + 4}}\xi \right) \right|^2 \, \left| \sin (n\pi \xi) \right|^2.
%\ee
which implies (see \cite{at} and \cite{Rebar} for more details), as in \cite{aht} for the only one string equation, 
that the system \eqref{Ap1b}-\eqref{Ap1c} is not exponentially stable in the energy space for all $\xi$ and $\beta.$
%But according to \rfb{obsineq1}, the main theorem of Arendt-Batty-Lyubich-Phong
%\cite{arentbatty, engel,TW}, the strong stabilization in the energy space is satisfied if and only if 
%$$\xi \notin \mathbb{Q}, \, \frac{\beta \xi}{\sqrt{\beta^2 + 4}} \notin \left\{\frac{2p - 1}{2q}, \, 
%p,q \in \mathbb{N}^* \right\},
%$$
%see  \cite{at1} and \cite{aht} for more details.

\subsection{Second example : mixed boundary conditions}
Consider the following  coupled wave equations
\begin{align}\label{Ap1bb}
&\ddot w_1(t,x)-\frac{\partial^2 w_1}{\partial x^2}(t,x) + w_1(t,x) + \alpha_1 \dot w_1(t,\xi)\delta_\xi+\alpha_2 \dot
 w_1(t-\tau,\xi)\delta_\xi+\beta \frac{\partial \dot w_2}{\partial x}(t,x)=0,\quad (t,x)\in(0,\infty)\times (0,1),\\
&\ddot w_2(t,x)-\frac{\partial^2 w_2}{\partial x^2}(t,x) + w_2(t,x) + \beta \frac{\partial \dot w_1}{\partial x}(t,x)=0,\quad
 (t,x)\in(0,\infty)\times (0,1),\\
&\frac{\partial w_1}{\partial x} (t,0) = \frac{\partial w_1}{\partial x}(t,1) = 0, \, w_2(t,0)=w_2(t,1)=0,\quad  t\in (0,\infty), \\
&w_i(0,x)=w_i^0(x),\,\,\dot w_i(0,x)=w_i^1(x),\,\,\dot w_1(s,x)=f_0(s,x), -\tau\leq s < 0,x\in (0,1), i=1,2.
\end{align}
with $\xi\in (0,1)$, $\beta$ is a positive constant and  $0<\alpha_2<\alpha_1$.

To put this control system into the framework of this paper, consider the spaces  $H_1=H_2=L^2(0,1)$ and the operators
$
A_1=A_2=-\frac{d^2}{dx^2}+I,
$
with  domains   $$\cD(A_1)= \left\{ u \in H^2(0,1), \, \frac{du}{dx}(0) = 0, \, \frac{du}{dx}(1) = 0 \right\},\, 
\cD(A_2)=H^2(0,1)\cap H^1_0(0,1)$$
  which are obviously  self-adjoint positive operators. In this case,  the domain of the  fractional power  
operators are given by
$$
\cD(A^\frac12_1)= H^1(0,1), \, \cD(A^\frac12_2)=H^1_0(0,1).
$$The operator $B$ and its adjoint $B^*$  are given by
$$
Bk= k\delta_\xi,\;\; k\in\mathbb{R}, \;\;B^*\varphi=\varphi(\xi), \;\; \varphi\in H^{1}(0,1)
$$ and  finally
$$
C=\beta\frac{d}{dx}, \mbox{ with } \cD(C)=H^1_0(0,1).
$$
It is clear that  $B^*:H^1(0,1)\rightarrow \mathbb{R}$ is  bounded and $C^*= - \beta\frac{d}{dx}$  with
\be H^1_0(0,1)\hookrightarrow \cD(C^*)= H^1(0,1).\ee
Assuming  $\beta<1$, as in the first example,  the condition \eqref{condition} is satisfied  with 
constant $\delta=\frac{\beta}{2}$.
Let us verify the boundedness of the transfer function of the above system. 
%We can easily verify that
%$$
%\left[\lambda^2I+A_1+\lambda^2C(\lambda^2+A_2)^{-1}C^*\right]^{-1}=
%\frac12\left[\lambda^2I+A_1+\lambda C\right]^{-1}
%+\frac12\left[\lambda^2I+A_2+\lambda C^*\right]^{-1}.
%$$
%Thus, we have the following decomposition of the transfer function
%$$
%H(\lambda)=\frac\lambda2B^*\left[\lambda^2I+A_1+\lambda C\right]^{-1}B+\frac\lambda2B^*
%\left[\lambda^2I+A_2+\lambda C^*\right]^{-1}B:= H_1(\lambda)+H_2(\lambda)
%$$
For this, let  $k\in \R$,  $Re \, \lambda > 0$ and the elleptic  system
\begin{align}
     & \lambda^2\phi_1(x)-\frac{d^2 \phi_1}{dx^2}(x) + \phi_1(x) + \lambda\beta\frac{d \phi_2}{dx}(x)= k \, \delta_\xi, \quad  x\in (0,1),
      \label{e1bv}\\
     & \lambda^2\phi_2(x)-\frac{d^2 \phi_2}{dx^2}(x) + \phi_2(x) + 
\lambda\beta\frac{d \phi_1}{dx}(x)=0,  \quad x\in (0,1),
     \label{e1bbv} \\
   & \frac{d\phi_1}{dx}(0)=\frac{d\phi_1}{dx}(1)=0, 
\label{e2bv} \\
   & \phi_2(0)=\phi_2(1)=0. \label{e1bbbv}
\end{align}
Then, 
$$
H(\lambda)= \lambda \, \phi_1(\xi) = \frac{\lambda}{2} \, \psi_1(\xi) +
 \frac{\lambda}{2} \, \psi_2(\xi):= H_1(\lambda) + H_2(\lambda) ,
$$
where $\psi_1 = \phi_1 - \phi_2, \, \psi_2 = \phi_1 + \phi_2$ satisfy the following equations
\begin{align}
     & \lambda^2\psi_1(x)-\frac{d^2 \psi_1}{dx^2}(x) + \psi_1 + \lambda\beta\frac{d \psi_1}{dx}(x)=0, \;\;
 x\in (0,\xi)\cup(\xi,1) \label{e1b}\\
     & \lambda^2\psi_2(x)-\frac{d^2 \psi_2}{dx^2}(x) + \psi_2 - \lambda\beta\frac{d \psi_2}{dx}(x)=0, \;\; x\in (0,\xi)\cup(\xi,1)
     \label{e1bb} \\
   & \frac{d(\psi_1 + \psi_2)}{dx}(0)=\frac{d(\psi_1 + \psi_2)}{dx}(1)=0, \label{e2b} \\
   & (\psi_1 - \psi_2)(0)=(\psi_1 - \psi_2)(1)=0, \label{e1bbb} \\
&[\psi_i]_\xi=0, \left[\frac{d \psi_i}{dx}\right]_\xi=k, \, i=1,2. \label{e3b}
\end{align}
Let $r_1, r_2$ be the roots of the equation $r^2-\beta\lambda r-\lambda^2 - 1 =0$,  which are
$
\frac{\beta \lambda}{2} \pm \sqrt{\frac{\beta^2 \lambda^2}{4} + \lambda^2 + 1}
$.
Then the solution of the equations \eqref{e1b}-\eqref{e1bbb} is given by
$$
\psi_1(x)=\left\{
  \begin{array}{ll}
    A_1 \, e^{r_1x} + B_1 \, e^{r_2x}, & x\in(0,\xi) \\
    C_1\,e^{r_1(x-1)} + D_1 \, e^{r_2(x-1)}, & x\in(\xi,1)
  \end{array}
\right.
$$
and
$$
\psi_2(x)=\left\{
  \begin{array}{ll}
    A_1\, e^{-r_1x} + B_1 \, e^{-r_2x}, & x\in(0,\xi) \\
    C_1 \,e^{-r_1(x-1)}+ D_1 \, e^{-r_2(x-1)}, & x\in(\xi,1).
  \end{array}
\right.
$$
Therefore,   \eqref{e3b} yields
$$
\psi_1(x)=\frac{k}{r_1-r_2}\left\{
  \begin{array}{ll}
  \frac{e^{-2r_1(\xi -1)} + 1}{e^{-r_1 \xi} - e^{r_1 (- \xi + 2)}} \, e^{r_1x} -
\frac{e^{-2r_2(\xi -1)} + 1}{e^{-r_2 \xi} - e^{r_2 (- \xi + 2)}} \, e^{r_2x}, & x \in (0,\xi),
\\

\left( e^{r_1} \, \frac{e^{-2r_1(\xi -1)} + 1}{e^{-r_1 \xi} - e^{r_1 (- \xi + 2)}} + e^{- r_1(\xi -1)} \right)\, e^{r_1(x-1)}  \\

-\left(e^{r_2} \, \frac{e^{-2r_2(\xi -1)} + 1}{e^{-r_2 \xi} - e^{r_2 (- \xi + 2)}} + e^{-r_2(\xi -1)} \right)\, e^{r_2(x-1)},
 & x\in(\xi,1)
  \end{array}
\right.
$$
and
$$
\psi_2(x)= \frac{k}{r_1 - r_2}\left\{
  \begin{array}{ll}
  \frac{e^{-2r_1(\xi -1)} + 1}{e^{-r_1 \xi} - e^{r_1 (- \xi + 2)}} \, e^{-r_1x} -
\frac{e^{-2r_2(\xi -1)} + 1}{e^{-r_2 \xi} - e^{r_2 (- \xi + 2)}} \, e^{-r_2x}, & x \in (0,\xi),
\\

\left( e^{r_1} \, \frac{e^{-2r_1(\xi -1)} + 1}{e^{-r_1 \xi} - e^{r_1 (- \xi + 2)}} + e^{- r_1(\xi -1)} \right)\, e^{-r_1(x-1)} \\
-\left( e^{r_2} \, \frac{e^{-2r_2(\xi -1)} + 1}{e^{-r_2 \xi} - e^{r_2 (- \xi + 2)}} + e^{-r_2(\xi -1)} \right)\, e^{-r_2(x-1)}, & x\in(\xi,1). \end{array}
\right.
$$
Consequently
$$
H_1(\lambda)= \frac{1}{2 \sqrt{\beta^2 + 4}} \, \frac{- \cosh[r_1(\xi - 1)] \, \sinh(r_2) \,  e^{r_1 \xi} + \cosh[r_2 (\xi -1)] \, \sinh(r_1) \, e^{r_2 \xi}}{\sinh(r_1) \, \sinh(r_2)},
$$
$$
H_2(\lambda)= \frac{1}{2 \sqrt{\beta^2 + 4}} \, \frac{- \cosh[r_1(\xi - 1)] \, \sinh(r_2) \, e^{-r_1 \xi} + \cosh[r_2(\xi - 1)] \, \sinh(r_1) \, e^{-r_2 \xi}}{\sinh(r_1) \, \sinh(r_2)}.
$$
As  $r_1$ and $r_2$ behave asymptotically as  $r_3: = \frac{\beta \lambda}{2} + \frac{\lambda}{2} \, 
\sqrt{\beta^2 + 4}$ and $r_4 := \frac{\beta \lambda}{2} - \frac{\lambda}{2} \, \sqrt{\beta^2 + 4}$, respectively
 it suffieses to see that for $r_3,r_4$, one has
$$
\sup_{Re\lambda=2\gamma}|H_1(\lambda)|\leq
\frac{1}{ \sqrt{\beta^2+4}}\, \frac{\cosh[\gamma(\xi - 1)(\beta + \sqrt{\beta^2 + 4})] \, \cosh[\gamma (\beta + \sqrt{\beta^2 + 4})] \, e^{\gamma \xi (\beta + \sqrt{\beta^2 +4})}}{\sinh[\gamma (\beta + \sqrt{\beta^2 + 4})] \, \sinh[\gamma (-\beta + \sqrt{\beta^2 + 4})]}.
$$
By similar calculus, we have the boundedness of $H_2$, and this achieves the claim.

Consider the conservative adjoint  system
\begin{align}\notag
&\frac{\partial^2 \phi}{\partial t^2}(t,x)-\frac{\partial^2 \phi}{\partial x^2}(t,x) + \phi(t,x) + 
\beta \frac{\partial^2 \psi}{\partial x \partial t}(t,x)=0,\quad (t,x)\in(0,\infty)\times (0,1),\\
&\frac{\partial^2 \psi}{\partial t^2}(t,x)-\frac{\partial^2 \psi}{\partial x^2}(t,x) + \psi(t,x) + \beta \frac{\partial^2 \phi}{\partial 
x \partial t}(t,x)=0,\quad (t,x)\in(0,\infty)\times (0,1),\notag\\
&\frac{\partial \phi}{\partial x}(t,0)=\psi(t,0)=\frac{\partial \phi}{\partial x}(t,1)=\psi(t,1)=0,\quad  
t\in (0,\infty), \notag\\
&\phi(0,x)=\phi^0(x),\,\,\dot \phi(0,x)=\phi^1(x),\,\, \psi(0,x)=\psi^0(x),\,\,\dot \psi(0,x)=\psi^1(x),\,\, x\in (0,1).
\notag
\end{align}
The  initial conditions can be written as 
\begin{align*}
\phi^0(x)=\sum_{n \in \mathbb{Z}^*} a_n \, \cos(n \pi x), \phi^1(x)=\sum_{n \in \mathbb{Z}^*} \lambda_n \, a_n \, \cos(n \pi x)\\
\psi^0(x)=\sum_{n \in \mathbb{Z}^*} a_n \, \sin(n \pi x), \psi^1(x)=\sum_{n \in \mathbb{Z}^*} \lambda_n \,  a_n \, \sin(n \pi x)
\end{align*}
with $\lambda_n = 
i \frac{2n\pi \beta}{\beta^2 + 4} \pm i \sqrt{\frac{4n^2 \pi^2 \beta^2}{(\beta^2 + 4)^2} + \frac{4 + 4 n^2 \pi^2}{\beta^2 + 4}},\, \, n \in \mathbb{Z}^*$,  and  $(\lambda_n \, a_n)\in l^2(\mathbb{C})$.
Hence, by standard technics, we obtain
$$
\phi(t,x) =\sum_{n \in \mathbb{Z}^*}
a_n \, e^{\lambda_n t} \,
\cos(n \pi x),
$$
and then 
\begin{align*}
\frac{\partial \phi}{\partial t}(t,\xi)=  \sum_{n \in \mathbb{Z}^*} \lambda_n \, a_n \, e^{\lambda_n t} \, \cos(n \pi \xi).
\end{align*}
Now, by the Ingham's inequality, for any  $T> \frac{\beta^2 + 4}{\beta + \sqrt{\beta^2+4}}$ there is  $C_{T,\xi,\beta}>0$ such that
\be
\label{obsineq1bb}
\int_0^T\left|\frac{\partial \phi}{\partial t}(t,\xi)\right|^2 dt
\geq C_{T,\xi,\beta} \sum_{n \in \mathbb{Z}^*} |\lambda_n|^2 \, |a_n|^2 \, \left| \cos (n \pi \xi) \right|^2.\notag
\ee
Finally, this  implies, as in \cite{abmb,aht} for the only one string equation, that the system is exponentially 
stable in the energy space  
if and only if $\xi$ is a rational number with coprime factorisation $\xi = \frac{p}{q}$, where $p$ is odd.

% ----------------------------------------------------------------
\bibliographystyle{amsplain}

\end{document}